\documentclass[10pt,leqno,a4paper]{amsart}

\usepackage{amsmath}
\usepackage{amssymb}
\usepackage{amsthm}
\usepackage{ifthen}
\usepackage{enumerate.mod}
\usepackage{calc}

\newtheorem{thm}{Theorem}
\newtheorem*{tthm}{\large Theorem}
\newtheorem{lem}[thm]{Lemma}

\newtheorem{cor}[thm]{Corollary}

\newtheorem*{tquestion}{Question}

\theoremstyle{definition}
\newtheorem{defn}[thm]{Definition}

\theoremstyle{remark}
\newtheorem{remark}[thm]{Remark}

\newcommand{\ch}{\mathrm{CH}}

\newcommand{\ma}{\textrm{\rm MA}}
\newcommand{\zfc}{\textrm{\rm ZFC}}

\newcommand{\pstar}{(*)}

\newcommand{\dom}{\operatorname{dom}}

\newcommand{\Lim}{\operatorname{Lim}}
\newcommand{\Fin}{\operatorname{Fin}}

\newcommand{\power}{\P}
\newcommand{\lbrak}{\bigl\|}
\newcommand{\rbrak}{\bigr\|}
\newcommand{\lmeas}{\mu\bigl(\lbrak}
\newcommand{\rmeas}{\rbrak\bigr)}

\newcommand{\otp}{\operatorname{otp}}

\renewcommand{\div}{\mathbin{/}}

\newcommand{\subseteqfnt}{\mathrel{\subseteq^*}}
\newcommand{\subsetfnt}{\mathrel{\subset^*}}

\newcommand{\divrel}[2]{\mathrel{{#1}\div{#2}}}

\newcommand{\forces}{\mskip 5mu plus 5mu\|\hspace{-2.5pt}{\textstyle \frac{\hspace{4.5pt}}{\hspace{4.5pt}}}\mskip 5mu plus 5mu}

\newcommand{\Iff}{\espc\mathrm{iff}\espc}

\renewcommand{\And}{\espc\text{and}\espc}

\newcommand{\add}{\operatorname{add}}

\newcommand{\cantorcube}[1]{\{0,1\}^{\hsp#1}}
\newcommand{\two}{\{0,1\}}

\newcommand{\pN}{\power(\N)}
\newcommand{\reals}{\mathbb R}

\newcommand{\N}{\mathbb N}
\newcommand{\integers}{\mathbb Z}

\newcommand{\nulls}{\mathcal N}

\newcommand{\random}{\R}

\DeclareFontFamily{U}{cmsy}{}
\DeclareFontShape{U}{cmsy}{m}{n}{<12> sfixed * [10] cmsy10 
<10> <9> <8> <7> <6> <5> sfixed * [10] cmsy10}{}
\DeclareSymbolFont{customtwo}{U}{cmsy}{m}{n} 
\DeclareMathSymbol{\sctn}{\mathord}{customtwo}{"78}

\DeclareFontFamily{U}{cmmi}{}
\DeclareFontShape{U}{cmmi}{m}{n}{<20> sfixed * [11] cmmib10 <12> sfixed * [10]
cmmi10 <10> <9> <8> sfixed * [6] cmmi6 <5> <6> <7> sfixed * [5] cmmi5}{}
\DeclareSymbolFont{custom}{U}{cmmi}{m}{n}
\DeclareMathSymbol{\rharpoon}{\mathord}{custom}{"2A}
\newlength{\widt}
\newlength{\widttwo}
\newlength{\hgt}

\newcommand{\oone}{{\omega_1}}
\newcommand{\vomega}{\varOmega}

\newcommand{\<}{\langle}
\renewcommand{\>}{\rangle}

\newcommand{\espc}{\quad}

\newcommand{\ulc}{\ulcorner}
\newcommand{\urc}{\urcorner}

\newcommand{\B}{\mathcal B}

\newcommand{\F}{\mathcal F}

\renewcommand{\L}{\mathcal L}
\newcommand{\G}{\mathcal G}

\renewcommand{\P}{\Pcal}
\newcommand{\Pcal}{\mathcal P}
\newcommand{\Q}{\mathcal Q}
\newcommand{\R}{\mathcal R}

\newcommand{\hsp}{\mspace{1.5mu}}

\newcounter{saveenumi}
\newcommand{\save}{\setcounter{saveenumi}{\value{enumi}}}
\newcommand{\restore}{\setcounter{enumi}{\value{saveenumi}}}

\newcommand{\pfn}{\dashrightarrow}

\newcommand{\Ult}{\operatorname{Ult}}

\newcommand{\xor}{\veebar}

\begin{document}

\title{On the strength of Hausdorff's gap condition}
\author{James Hirschorn}
\date{September 21, 2002}
\address{Institut f\"ur Formale Logik, Universit\"at Wien, Austria}
\email{James.Hirschorn@logic.univie.ac.at}
\thanks{Research supported by Lise Meitner Fellowship, 
Fonds zur F\"orderung der wissenschaftlichen Forschung (FWF)}
\urladdr{http://www.logic.univie.ac.at/\textasciitilde hirschor/}

\begin{abstract}
Hausdorff's gap condition was satisfied by his original 1936 construction of an
$(\oone,\oone)$ gap in $\pN\div\Fin$. We solve an open problem in determining
whether Hausdorff's condition is actually stronger than the more modern
indestructibility condition, by constructing an indestructible $(\oone,\oone)$ gap
not equivalent to any gap satisfying Hausdorff's condition,
from uncountably many random reals.
\end{abstract}

\maketitle

\section{Introduction}
\label{sec:introduction}

%\the\textwidth\the\linewidth

A \emph{pregap} in a Boolean algebra $(\B,\le)$ is an orthogonal pair $(A,B)$ of
subsets of $\B$, i.e.~%
\begin{enumerate}[(i)]
\item $a\cdot b=0$ for all $a\in A$ and $b\in B$,
\save
\end{enumerate}
and it is a \emph{gap} if additionally there is no element $c$ of $\B$ such that
\begin{enumerate}[(i)]
\restore
\item $a< c$ for all $a\in A$, and $b< -c$ for all $c\in B$.
\save
\end{enumerate}
Such an element $c$ is said to \emph{interpolate} the pregap.
A \emph{linear pregap} is a pregap $(A,B)$ where both $A$ and $B$ are linearly
ordered by $\le$, and for a pair of linear order types $(\varphi,\psi)$, a
$(\varphi,\psi)$ \emph{pregap} in a Boolean algebra $(\B,\le)$ is a linear pregap
$(A,B)$ where $\otp(A,\le)=\varphi$ and $\otp(B,\le)=\psi$. Thus $(A,B)$ is a
$(\varphi,\psi)$ \emph{gap} if it is a $(\varphi,\psi)$ pregap for which no element of
$\B$ can be used to extend $(A,B)$ to a $(\varphi+1,\psi)$ pregap or a
$(\varphi,\psi+1)$ pregap.

Gaps in the Boolean algebra $(\pN\div\Fin,\subseteqfnt)$ have a long history with
some basic results appearing as early as 1873, including Hadamard's theorem~\cite{Had} that
there are no $(\delta,\delta)$ gaps in $\pN\div\Fin$ for any ordinal $\delta$ with
countable cofinality. Indeed, one of the major achievements in early Set Theory was
Hausdorff's construction~\cite{Ha} in 1936 of an $(\oone,\oone)$ gap. 

While being a pregap in $\pN\div\Fin$ is absolute for transitive models, the
property of being a gap is not. For example, if $(A,B)$ is an $(\oone,\oone)$ gap in
$\pN\div\Fin$ and $\Q$ is a poset which collapses $\oone$, then by Hadamard's
theorem forcing with $\Q$ must introduce an element of $\pN\div\Fin$ which
interpolates $(A,B)$ and thus renders it a non-gap. Avoiding this particular
example, an $(\oone,\oone)$ pregap $(A,B)$ in $\pN\div\Fin$ is called
\emph{destructible} if there is an $\oone$-preserving poset which forces that
$(A,B)$ is not a gap. 

Destructibility is in fact a Ramsey theoretic phenomenon.
This becomes clear when one considers the characterization
of destructibility below (Lem\-ma~\ref{l-2}).  
When working with pregaps in
$\pN\div\Fin$ one often works with representatives (i.e.~subsets of $\N$) of the
equivalence classes. In this case, whenever $a_\alpha,b_\alpha\subseteq\N$ ($\alpha<\delta$) is a
$(\delta,\delta)$ pregap in $\pN\div\Fin$ for some ordinal $\delta$, we may
assume---in order to avoid trivialities and thereby obtain more concise
results---that the representatives have been chosen so that
\begin{equation}
  \label{eq:8}
  a_\alpha\cap b_\alpha=\emptyset\espc\text{for all $\alpha<\delta$}.
\end{equation}
We also assume that the enumerations respect the well orderings of both
$\{a_\alpha:\alpha<\delta\}$ and $\{b_\alpha:\alpha<\delta\}$ by $\subsetfnt$.

Given a pregap $(a_\alpha,b_\alpha:\alpha<\oone)$, define a partition of
$[\oone]^2=K_0\cup K_1$ by 
\begin{equation}
  \label{eq:7}
  \{\alpha,\beta\}\in K_0\Iff (a_\alpha\cap b_\beta)\cup(a_\beta\cap b_\alpha)=\emptyset.
\end{equation}

\begin{lem}
\label{l-2}
For every $(\oone,\oone)$ pregap in $\pN\div\Fin$ the following
are equivalent\textup:
\textup{\begin{enumerate}[(a)]
\item \textit{$(a_\alpha,b_\alpha:\alpha<\oone)$ is destructible.}
\item \textit{There is no uncountable $1$-homogeneous subset of $\oone$.}
\item \textit{There exists a poset with the ccc forcing that
  $(a_\alpha,b_\alpha:\alpha<\oone)$ is not a gap.}
\end{enumerate}}
\end{lem}
\begin{proof}
See~\cite{Sch}.
\end{proof}

The existence of destructible $(\oone,\oone)$ gaps in $\pN\div\Fin$ is independent
of the usual axioms of mathematics ($\zfc$). For example, it is a theorem of Kunen~\cite{Ku}
that $\ma_{\aleph_1}$ implies that all $(\oone,\oone)$ gaps are indestructible,
while theorems of Todor\v cevi\'c are that a destructible $(\oone,\oone)$ gap can
be constructed from either a diamond sequence or a Cohen real
(see~\cite{D},~\cite{TF}, resp.).

The existence of indestructible $(\oone,\oone)$ gaps 
can be proved from the usual axioms of mathematics. Indeed
the gap originally constructed by Hausdorff was indestructible, and in fact
satisfied an even stronger condition:
\begin{equation}
  \label{eq:36}
  \{\alpha<\beta:a_\alpha\cap b_\beta\setminus k=\emptyset\}\text{ is
    finite}\espc\text{for all $k\in\N$, for all $\beta<\oone$.}
\end{equation}

Hausdorff's condition is strictly stronger than indestructibility because one can
always modify the initial parts $(a_\alpha,b_\alpha:\alpha<\omega+1)$ of an
indestructible gap to make~\eqref{eq:36} fail. However, the following definition
is clearly a correct description of the equivalence of two gap structures.

\begin{defn}
Two gaps $(A,B)$ and $(A',B')$ in some Boolean algebra $(\B,\le)$ are
\emph{equivalent} if $A$ is $\le$-cofinal in $A'$, $A'$ is $\le$-cofinal in $A$, $B$
is $\le$-cofinal in $B'$ and $B'$ is $\le$-cofinal in $B$.
\end{defn}

\noindent For example, if $(A,B)$ is an $(\oone,\oone)$ pregap in some Boolean
algebra then for any uncountable $A'\subseteq A$ and $B'\subseteq B$, 
$(A',B')$ is equivalent to the original pregap. 

Hausdorff's condition is not just of
historical interest. In some contexts it is a more natural condition than
indestructibility. For example, as shown in~\cite{rst:A}, $\pstar$---a combinatorial
principle for $P$-ideals of countable sets which is compatible with $\ch$---implies that all
$(\oone,\oone)$ gaps in $\pN\div\Fin$ are indestructible because for a given
$(\oone,\oone)$ gap $(a_\alpha,b_\alpha:\alpha<\oone)$ the ideal of all
countable subsets $\vomega$ of $\oone$ on which
$(a_\alpha,b_\alpha:\alpha\in\vomega)$ satisfies Hausdorff's condition forms a $P$-ideal, and
applying $\pstar$ to this $P$-ideal establishes that $(A,B)$ is equivalent with a
gap satisfying Hausdorff's condition. 

Scheepers in~\cite{Sch} and~\cite{Sch1} has asked whether indestructibility is as
strong as Hausdorff's condition. 

\begin{tquestion}[Scheepers]
Is every indestructible $(\oone,\oone)$ gap in $\pN\div\Fin$ equivalent with an
$(\oone,\oone)$ gap satisfying Hausdorff's condition\textup?
\end{tquestion}

We give a negative answer (i.e.~we show that the statement in Question is
consistently false).

\begin{tthm}
If $\random$ is a nonseparable measurable algebra, then with positive probability,
there exists an indestructible $(\oone,\oone)$ gap in
$(\pN\div\Fin,\subseteqfnt)$ which is not equivalent to any $(\oone,\oone)$ gap in
$\pN\div\Fin$ satisfying Hausdorff's condition.
\end{tthm}

This can be put in a different setting. It is a corollary that the classical
hypothesis that the Lebesgue measure can be extended to all subsets of $\reals$
distinguishes between indestructibility and Hausdorff's condition.
This is an immediate consequence of the
Theorem and known absoluteness results for forcing extensions by a
large enough measure algebra, from a real-valued measurable cardinal of size at most
continuum (see~\cite{TF}). 

\begin{cor}
If the Lebesgue measure on the real line can be extended to a measure whose domain
is all of $\P(\reals)$, then there exists a indestructible $(\oone,\oone)$ gap in
$\pN\div\Fin$ which is not equivalent to any gap satisfying Hausdorff's condition.
\end{cor}
\begin{proof}
Suppose that $\mu:\power(\reals)\to[0,\infty]$ is a measure extending the Lebesgue
measure, and let $\nulls_\mu$ be the ideal of all subsets $N\subseteq\reals$ with
$\mu(N)=0$. Let $\G$ be a generic filter over $V$ of the poset
$(\Q,\subseteq)$ where $\Q=\power(\reals)\setminus\nulls_\mu$. Since
\begin{equation}
  \label{eq:43}
  (V,\in,\G)\models\ulc \G\text{ is a $\sigma$-complete ultrafilter on $\reals$}\urc,
\end{equation}
in $V[\G]$, the ultrapower $\Ult(V,\G)=(V^{\reals\cap V}\cap V)\div\G$ ($f,g\in
V^{\reals\cap V}\cap
V$ are equivalent if $\{r\in\reals\cap V:f(r)=g(r)\}\in\G$) with its ordering ${\in^\reals}\div\G$
(i.e.~$\{r:f(r)\in g(r)\}\in\G$) is well-founded, and thus there is an
elementary embedding 
\begin{equation}
  \label{eq:65}
  j:V\to M\cong \Ult(V,\G)
\end{equation}
where $M$ is the transitive (Mostowski) collapse of
the ultrapower. 

Arguing in $V$: Note that for every $\Q$-name $\dot f$ for a
member of $(V^\reals)\check{\vphantom)}$ there is a $g_{\dot f}\in V^\reals$ such that 
$\Q\forces[\dot f]=[\check g_{\dot f}]$ (where $[f]$ denotes the
equivalence class of $f$ in $\Ult(V,\G)$). Suppose that $\{\dot f:f\in\F\}$ is a
family of cardinality less than $\add(\nulls_\mu)$ of $\Q$-names for
members of $(V^\reals)\check{\vphantom)}$.  Defining $G\in V^\reals$ by 
$G(r)=\{g_{\dot f}(r):\dot f\in\F\}$, for every $\dot f\in\F$, $\Q\forces[\dot
f]\divrel{\in^{\reals}}{\dot\G}[\check G]$ since
$\{r\in\reals:g_{\dot f}(r)\notin G(r)\}=\emptyset\in\nulls_\mu$; 
and conversely, if $A\in\Q$ forces that
$[\dot h]\ne[\dot f]$ for all $\dot f\in\F$, then $\mu(\{r\in A:g_{\dot h}(r)=g_{\dot f}(r)\})=0$
for all $\dot f\in\F$, and thus $\mu(\{r\in A:g_{\dot h}(r)\in G(r)\})=0$ since
$|\F|<\add(\nulls_\mu)$, and therefore $A\forces[\dot
h]\divrel{\notin^\reals}{\dot\G}[\check G]$. This proves that
\begin{equation}
  \label{eq:64}
  [M]^{<\add(\nulls_\mu)^V}\subseteq M.
\end{equation}

Forcing with $\Q$ is the same as forcing with $\random^+$ where $(\random,\nu)$ is
the measure algebra of the measure space $(\reals,\P(\reals),\mu)$,
i.e.~$\random=\power(\reals)\div\nulls_\mu$ and 
$A\in\power(\reals)\setminus\nulls_\mu\mapsto [A]\in\random^+$
is a surjective embedding. By considering if necessary a homogeneous principle ideal
$\random_z=\{x\in\random:x\le z\}$ $(z\in\random^+)$ of $\random$, Maharam's theorem states that
the probability algebra $(\random_z,\mu_z)$, where $\mu_z(x)=\mu(x)\div\mu(z)$ for
all $x\in\random_z$, is isomorphic (as a measure algebra) 
to the measure algebra of $\cantorcube\theta$ with its Haar
probability measure, where $\theta$ is either 0 or an infinite cardinal and is called the
\emph{Maharam type} of $\random_z$ (see~\cite{F3}). In other
words, $\random_z^+$ is isomorphic to the canonical poset for adding $\theta$ random reals. It is a
theorem of Gitik and Shelah~\cite{GS} that the measure algebra of a measure space with domain
$\P(\reals)$ has Maharam type $2^{2^{\aleph_0}}$. In particular, $\random$ is
nonseparable, and therefore by the Theorem, there is an $A\in\Q$ forcing
the existence of an indestructible $(\oone,\oone)$ gap $(A,B)$ in $\pN\div\Fin$
which is not equivalent to any gap satisfying Hausdorff's condition. Assume
that $A\in\G$.

It is a classical theorem of Ulam that $\add(\nulls_\mu)\le2^{\aleph_0}$ is at least as large as
the smallest weakly inaccessible cardinal, and in particular~\eqref{eq:64} implies that, in $V[\G]$,
$H_{\aleph_2}\subseteq M$. Since $(\oone,\oone)$ pregaps in $\pN\div\Fin$ are
objects of $H_{\aleph_2}$, by elementarity there exists in $V$ such a gap.
\end{proof}

\begin{remark}
It follows from some of the proofs in~\cite{HI} that, assuming $\ma_{\aleph_1}$,
in the extension by a separable measurable algebra, every $(\oone,\oone)$ gap in
$\pN\div\Fin$ is equivalent to a Hausdorff gap. Thus the construction is not
possible with only one random real.
\end{remark}

\section{Construction}
\label{sec:construction}

\subsection{Measure theoretic characterization}
\label{sec:meas-theor-char}

We have the following simple necessary condition for the Hausdorff property of an
$\random$-name for pregap.

\begin{lem}\label{l-1}
Let $(\random,\mu)$ be a probability algebra and
let $(\dot a_\alpha,\dot b_\alpha:\alpha<\oone)$ be an $\random$-name for an
$(\oone,\oone)$ pregap in $\pN\div\Fin$.
If for some $h\in c_0$, 
\begin{equation}
  \label{eq:1}
  \lmeas\dot a_\alpha\cap\dot b_\beta\setminus k\ne\emptyset\rmeas\le h(k)
\end{equation}
for all $\alpha,\beta<\oone$ and all $k\in\N$, then with probability one, $(\dot
a_\alpha,\dot b_\alpha:\alpha<\oone)$ is not equivalent to any gap
satisfying Hausdorff's condition.
\end{lem}
\begin{proof}
First note that if $(a_\alpha,b_\alpha:\alpha<\oone)$ is an $(\oone,\oone)$ pregap
in $\pN\div\Fin$, and $(a'_\alpha,b'_\alpha:\alpha<\oone)$ is an equivalent pregap,
then there must be an uncountable $X\subseteq\oone$, $\alpha_\xi\ge\xi$ for all
$\xi\in X$ and an integer $l$ such that $a'_\xi\setminus l\subseteq a_{\alpha_\xi}$
and $b'_\xi\setminus l\subseteq b_{\alpha_\xi}$ for all $\xi\in X$. Let
$X_0\subseteq X$ be uncountable where the sequence $\alpha_\xi$ $(\xi\in X)$ is
strictly increasing
on $X_0$. Then if $(a'_\alpha,b'_\alpha:\alpha<\oone)$ satisfies Hausdorff's
condition, so does $(a_{\alpha_\xi},b_{\alpha_\xi}:\xi\in X_0)$. 

Suppose then that $\dot Y$ is an $\random$-name for an uncountable subset of
$\oone$, and that $x\in\random^+$ forces that $(\dot a_\alpha,\dot
b_\alpha:\alpha\in \dot Y)$ has Hausdorff's property. 
Then there is a $\delta>0$ and an uncountable $X\subseteq \oone$ such that
\begin{equation}
  \label{eq:9}
  \mu\bigl(x\cdot\lbrak\alpha\in\dot Y\rmeas>\delta\espc\text{for all $\alpha\in X$}.
\end{equation}
Thus by Gillis' Theorem~\cite{G} there is an uncountable $X_0\subseteq X$ such that 
\begin{equation}
  \label{eq:11}
  \mu\bigl(x\cdot\lbrak\alpha\in\dot Y\rbrak\cdot\lbrak\beta\in \dot Y\rmeas>\delta^2
   \espc\text{for all $\alpha,\beta\in X_0$}.
\end{equation}
Let $k$ be large enough so that $h(k)\le\delta^2\div2$.
Pick $\beta\in X_0$ so that $X_0\cap\beta$ is infinite. Since
$x\cdot\lbrak\beta\in\dot Y\rbrak$ forces that there
are only finitely many $\alpha\in\dot Y\cap\beta$ such that $\dot a_\alpha\cap\dot
b_\beta\setminus k=\emptyset$, there is an $\bar\alpha\in X_0\cap\beta$ such that 
\begin{equation}
  \label{eq:12}
  \mu\bigl(x\cdot\lbrak\bar\alpha\in\dot Y\rbrak\cdot\lbrak\beta\in\dot Y\rbrak
  \cdot\lbrak\dot a_{\bar\alpha}\cap\dot b_\beta\setminus k=\emptyset\rmeas<\frac{\delta^2}2.
\end{equation}
Combining~\eqref{eq:11} and~\eqref{eq:12} yields $\lmeas\dot a_{\bar\alpha}\cap\dot
b_\beta\setminus k\ne\emptyset\rmeas>\delta^2\div2\ge h(k)$.
\end{proof}

\subsection{Notation}
\label{sec:notation}

We denote the set of all finite partial functions from $X$ into $Y$ by $\Fin(X,Y)$.
\emph{Logical and} (digitwise base 2) is denoted by `$\land$' and \emph{exclusive or} is denoted
by `$\xor$' (e.g.~for $i,j\in\{0,1\}$, $i\xor j=0$ if $i=j$, and $i\xor j=1$ if
$i\ne j$). Define a mapping $w:\oone\times\N\times\integers\to\Fin(\oone,\two)$ 
where $\dom(w(\alpha,i,j))=[\alpha,\alpha+i)$ and 
\begin{equation}
  \label{eq:4}
  w(\alpha,i,j)(\alpha+k)=\frac{j\mod{2^i}}{2^k}\land 1\espc\text{for all $k<i$,}
\end{equation}
in other words, the concatenation
\begin{equation}
w(\alpha,i,j)(\alpha+i-1)w(\alpha,i,j)(\alpha+i-2)\cdots w(\alpha,i,j)(\alpha)
\end{equation}
is the base $2$ representation of $j\mod{2^i}$. 

For a set $X$, we write $\bigl(\random_{(X)},\mu_{(X)}\bigr)$ for the measure algebra of the space
$\cantorcube X$ with its Haar probability measure. And for a finite partial function
$s:X\pfn\two$, we let $[s]_{(X)}\in\random_{(X)}$ be the equivalence class of the measurable
set
\begin{equation*}
  \bigl\{z\in\cantorcube X:z\supseteq s\bigr\}.
\end{equation*}
Sets of this form are called \emph{basic elements} of $\random_{(X)}$. It will be
convenient to be able to denote $0\in\random_{(X)}$ with the $[s]_{(X)}$ notation, and thus we adopt the
extension where 
\begin{equation}
  \label{eq:2}
  [\<-1\>]_{(X)}=0.
\end{equation}

\subsection{Proof of Theorem}
\label{sec:proof}

Define $g,h:\N\to\N$ by $g(n)=0\cdot 2^0+1\cdot 2^1+\cdots+n\cdot 2^n$, and 
\begin{equation}
  \label{eq:6}
  h(m)=\min_{n\in\N}m<g(n).
\end{equation}
The Theorem is proved by constructing an $\random_{(\oone)}$-name
$(\dot a_\alpha,\dot b_\alpha:\alpha<\oone)$ for a pregap such that 
\begin{alignat}{2}
  \label{eq:3}
  \lbrak \dot a_\alpha\cap\dot b_\alpha=\emptyset\rbrak&=1\espc
  &&\text{for all $\alpha$},\\
  \label{eq:5}
  \lbrak\dot a_\alpha\cap\dot b_\beta\ne\emptyset\rbrak&=1\espc
  &&\text{for all $\alpha\ne\beta$},\\
 \label{eq:50} 
  \mu_{(\oone)}\bigl(\lbrak m\in \dot a_\alpha\cap\dot b_\beta\rmeas&\le2^{-2h(m)}\espc
  &&\text{for all $m$, for all $\alpha,\beta$.}
\end{alignat}
Then by conditions~\eqref{eq:3} and~\eqref{eq:5}
and Lemma~\ref{l-2},  $(\dot a_\alpha,\dot b_\alpha:\alpha<\oone)$ is
an indestructible $(\oone,\oone)$ gap with probability one, and since $m\mapsto2^{-2h(m)}\in
\ell_1$, by condition~\eqref{eq:50} and applying Lemma~\ref{l-1} to
$k\mapsto\sum_{m=k}^\infty 2^{-2h(m)}\in c_0$, 
with probability one it is not equivalent to a
Hausdorff gap.

Suppose we are given a measurable algebra $\random$ with probability measure
$\mu$. Assuming $\random$ is nonseparable, i.e.~the measure algebra topology on
$\random$ induced by $\mu$ is nonseparable,
by Maharam's Theorem there is a measure algebra embedding of
$\bigl(\random_{(\oone)},\mu_{(\oone)}\bigr)$ into some principal ideal $(\random_z,\mu_z)$
($z\in\random^+$, $\mu_z(x)=\mu(x)\div\mu(z)$ for $x\in\random_z$) of
$\random$. Thus, forcing over the poset $\random^+$, $z$ forces that there exists an
indestructible $(\oone,\oone)$ gap not equivalent to any Hausdorff gap, 
completing the proof of the Theorem.

Since  an $\random_{(\oone)}$-name for an $(\oone,\oone)$ pregap
satisfying~\eqref{eq:3},~\eqref{eq:5} and~\eqref{eq:50} can clearly be expressed in the language
$\L_\omega(Q)$, we may use forcing with an $\oone$-preserving poset to perform
the construction. I.e.~if the existence of such an object is consistent then by
Keisler's Completeness Theorem for $\L_\omega(Q)$~\cite{Ke} it does exist. Indeed we are not
aware of how to do the construction differently, e.g.~by recursion on $\alpha<\oone$.

The following observation is used to satisfy~\eqref{eq:5}:
For every $n\ge1$, whenever $[\alpha,\alpha+n)\cap[\beta,\beta+n)=\emptyset$,
\begin{equation}
  \label{eq:13}
  \sum_{i=1}^n\sum_{j=0}^{2^i-1}[w(\alpha,i,j)\cup w(\beta,i,j\xor2^{i-1})]
  +\sum_{j=0}^{2^n-1}[w(\alpha,n,j)\cup w(\beta,n,j)]=1.
\end{equation}

Define $\<\>:(\N\setminus\{0\})\times\N\times\N\to\N$ by
\begin{equation}
  \label{eq:14}
  \<i,j,k\>=g(i-1)+ji+k.
\end{equation}
Note that $h(\<i,j,k\>)=i$ for all $i=1,2,\dots$, for all $j<2^i$, for all
$k<i$. For the remainder we write $\mu$ for $\mu_{(\oone)}$.

Let $\Q$ be the poset of all conditions
$(\vomega_p,\{s_p^\alpha,t_p^\alpha:\alpha\in\vomega_p\})$ where
\begin{enumerate}[(a)]
\item\label{item:1} $\vomega_p\subseteq\Lim(\oone)$ is finite,
\item\label{item:2} for some $n_p\in\N$, for all $\alpha\in\vomega_p$,
  \begin{equation*}
    s_p^\alpha,t_p^\alpha:g(n_p)\to\{\<-1\>\}\cup
    \bigcup_{\substack{\gamma\in\Lim(\oone)\\ \gamma\le\alpha}}
    \Fin\bigl([\gamma,\gamma+\omega),\{0,1\}\bigr),
  \end{equation*}
\item\label{item:6} $n_p\ge|\vomega_p|$,
\item\label{item:3} for all $i=1,\dots,n_p$, for all $j<2^{i}$, for all $k<i$,
  $\mu\bigl(\bigl[s_p^\alpha(\<i,j,k\>)\bigr]\bigr)\le 2^{-i}$ and 
  $\mu\bigl(\bigl[t_p^\alpha(\<i,j,k\>)\bigr]\bigr)\le2^{-i}$,
\item\label{item:4} for all $\alpha\in\vomega_p$, and all $i,j$, there exists
$k<i$ such that $s_p^\alpha(\<i,j,k\>)=w(\alpha,i,j)$ and
  $t_p^\alpha(\<i,j,k\>)=w(\alpha,i,j\xor2^{i-1})$, 
\item\label{item:7} for all $m<g(n_p)$, $[s_p^\alpha(m)]\cdot[t_p^\alpha(m)]=0$,
\item \label{item:9} for all $\alpha,\beta$ in $\vomega_p$, for all $i,j,k$,
  $\mu\bigl(\bigl[s^\alpha_p(\<i,j,k\>)\bigr]\cdot\bigl[t^\beta_p(\<i,j,k\>)\bigr]\bigr)\le2^{-2i}$,
\item\label{item:8} for all $\alpha\ne\beta$ in $\vomega_p$, there exists
  $E_{\alpha\beta}\subseteq g(n_p)$ such that 
\begin{enumerate}[(i)]
\item\label{item:10} $s_p^\alpha(m)\in\Fin\bigl([\alpha,\alpha+\omega),\two\bigr)$ 
  for all $m\in E_{\alpha\beta}$,
\item \label{item:11} $t_p^\beta(m)\in\Fin\bigl([\beta,\beta+\omega),\two\bigr)$ for all $m\in E_{\alpha\beta}$,
\item\label{item:12}   $\sum_{m\in E_{\alpha\beta}}[s_p^\alpha(m)]\cdot[t_p^\beta(m)]=1$,
\end{enumerate}
\save
\end{enumerate}
ordered by $q\le p$ if
\begin{enumerate}[(a)]
\restore
\item $\vomega_q\supseteq\vomega_p$,
\item for all $\alpha\in\vomega_p$, $s_q^\alpha\supseteq s_p^\alpha$ and
  $t_q^\alpha\supseteq t_p^\alpha$,
\item\label{item:5} for all $\alpha<\beta$ in $\vomega_p$, for all $m\in g(n_q)\setminus g(n_p)$,
  $[s^\alpha_q(m)]\le[s^\beta_q(m)]$ and
  $[t^\alpha_q(m)]\le[t^\beta_q(m)]$.
\end{enumerate}

\begin{lem}\label{l-3}
$\Q$ has precaliber $\aleph_1$ \textup(and in particular has the ccc\textup; in particular, $\Q$
does not collapse $\oone$\textup{).}
\end{lem}
\begin{proof}
Let $p(\xi)$ $(\xi<\oone)$ be an uncountable sequence of conditions. Let $X\subseteq\oone$ be an
uncountable set such that $\{\vomega_{p(\xi)}:\xi\in X\}$ forms a $\Delta$-system,
say with root $\vomega<\vomega_{p(\xi)}\setminus\vomega$ of size $l$, 
where $|\vomega_{p(\xi)}\setminus\vomega|=l^*$ for all
$\xi\in X$. By going to an uncountable subsequence, we may also assume that for all
$\xi\in X$,
\begin{alignat}{2}
  n_{p(\xi)}&=n,&&\\
  \label{eq:30}
  s_{p(\xi)}^\alpha&=s^\alpha\espc&&\text{for all $\alpha\in\vomega$,}\\
  \label{eq:31}
  t_{p(\xi)}^\alpha&=t^\alpha&&\text{for all $\alpha\in\vomega$.}
\end{alignat}

For each $\xi$, let $\{\gamma(\xi,0),\dots,\gamma(\xi,{l+l^*-1})\}$ be the
increasing enumeration of $\vomega_{p(\xi)}$.
By going to an uncountable subsequence we assume that 
\begin{alignat}{2}
  \label{eq:16}
  k_{\gamma(\xi,d)}(i,j)&=l_d(i,j)\espc&&\text{for all $d=0,\dots,l+l^*-1$, for all $\xi\in X$,}\\
  \label{eq:17}
  E_{\gamma(\xi,d)\gamma(\xi,{\bar d})}&=D_{d\bar d}&&\text{for all $d\ne\bar d$,
    for all $\xi\in X$},
\end{alignat}
where $k_{\gamma(\xi,d)}(i,j)<i$ satisfies the requirement of~\eqref{item:4}, for
$i=1,\dots,n$ and $j<2^i$, for the condition $p(\xi)$, and the sets 
$E_{\gamma(\xi,d)\gamma(\xi,{\bar d})}\subseteq g(n)$ satisfy the
requirements of~\eqref{item:8} for the condition $p(\xi)$.

For each $\alpha,\beta\in\oone$, let $\varphi_{\alpha\beta}:\oone\to\oone$ be the
bijection which swaps $\alpha+n$ with $\beta+n$ for all $n\in\N$, and fixes all
other ordinals. 
For $s\in\bigcup_{\alpha\in\Lim(\oone)}\Fin\bigl([\alpha,\alpha+\omega),\two\bigr)$,
define $\delta(s)\in\Lim(\oone)$ so that
$\dom(s)\subseteq[\delta(s),\delta(s)+\omega)$. By going to an uncountable
subsequence we can assume that there are sets $F_d\subseteq g(n)$ where
\begin{align}
  \label{eq:22}
  \bigl\{m<g(n):\delta\bigl(s^{\gamma(\xi,d)}_{p(\xi)}(m)\bigr)=\gamma(\xi,d)\bigr\}&=F_d,
\end{align}
for all $d=0,\dots,l+l^*$, for all $\xi\in X$; and that
\begin{alignat}{2}
  \label{eq:25}
  s^{\gamma(\xi,d)}_{p(\xi)}(m)\circ\varphi_{\gamma(\xi,d)\gamma(\eta,d)}&=s^{\gamma(\eta,d)}_{p(\eta)}(m)
  \espc&&\text{for all $m\in F_d$},
\end{alignat}
for all $d=0,\dots,l+l^*$, for all $\xi,\eta\in X$.

For each $\xi$, let 
$\varGamma_\xi=\bigl\{\delta\bigl(s_{p(\xi)}^\alpha(m)\bigr):\alpha\in\vomega_{p(\xi)}$, $m<g(n)\bigr\}$ and
$\varLambda_\xi=\bigl\{\delta\bigl(t_{p(\xi)}^\alpha(m)\bigr):\alpha\in\vomega_{p(\xi)}$,
$m<g(n)\bigr\}$. By going to an
uncountable subsequence we can assume that $\{\varGamma_\xi:\xi\in X\}$ and
$\{\varLambda_\xi:\xi\in X\}$ form $\Delta$-systems with roots $\varGamma$ and
$\varLambda$, respectively, and that
\begin{equation}
  \label{eq:23}
  \varGamma_\xi\cap\varLambda_\eta=\varGamma\cap\varLambda\espc\text{for all $\xi\ne\eta$.}
\end{equation}
Refining again, we can assume that there are sets $G_d\subseteq g(n)$ where 
\begin{equation}
  \label{eq:34}
  \bigl\{m<g(n):\delta\bigl(s_{p(\xi)}^{\gamma(\xi,d)}(m)\bigr)\in\varGamma\bigr\}=G_d
\end{equation}
for all $d$ and all $\xi$, and that
\begin{alignat}{2}
  \label{eq:24}
  s^{\gamma(\xi,d)}_{p(\xi)}(m)&=u(d,m)\espc&&\text{for all $m\in G_d$,}
\end{alignat}
for all $d$, for all $\xi$.

It remains to show that $\{p(\xi):\xi\in X\}$ is centered, but we can simplify
things by just proving that it is linked, because we only need the fact that $\Q$
does not collapse $\oone$. Fix $\xi\ne\eta$ in $X$. Put $n_q=\max(l+2l^*,n+1)$. 

By~\eqref{eq:30} and~\eqref{eq:31} we can define
$s_q^{\alpha},t_q^{\alpha}:g(n_q)\to\{\<-1\>\}\cup\bigcup_{\gamma\in\Lim(\alpha+1)}
\Fin([\gamma,\gamma+\omega),\two)$ so that 
\begin{alignat}{4}
  \label{eq:19}
  s_q^{\alpha}\restriction g(n)&=s_{p(\xi)}^{\alpha}&&\And t_q^{\alpha}\restriction g(n)&&=t_{p(\xi)}^\alpha
  \espc&&\text{for all $\alpha\in\vomega_{p(\xi)}$,}\\
  \label{eq:20}
  s_q^\alpha\restriction g(n)&=s_{p(\eta)}^\alpha&&\And
  t_q^\alpha\restriction g(n)&&=t_{p(\eta)}^\alpha
  &&\text{for all $\alpha\in\vomega_{p(\eta)}$;}
\end{alignat}
for all $i=n+1,\dots,n_q-1$, for all $j<2^{i}$, for all $k<i$: for all $d=0,\dots,l-1$,
\begin{alignat}{2}
  \label{eq:26}
  s_q^{\gamma(\xi,d)}(\<i,j,k\>)&=
  \left\{
  \begin{aligned}
    &w(\gamma(\xi,k),i,j),\\
    &\<-1\>,
  \end{aligned}
  \right.
\espc&&\begin{aligned}
  &\text{if $k\le d$,}\\
  &\text{if $k>d$,}
  \end{aligned}\\
  \label{eq:27}
  t_q^{\gamma(\xi,d)}(\<i,j,k\>)&=
  \left\{
  \begin{aligned}
    &w(\gamma(\xi,k),i,j\xor2^{i-1}),\\
    &\<-1\>,
  \end{aligned}
  \right.
\espc&&\begin{aligned}
  &\text{if $k\le d$,}\\
  &\text{if $k>d$,}
  \end{aligned}\\
\intertext{and for all $d=l,\dots,l+l^*-1$,}
  \label{eq:28}
  s_q^{\gamma(\xi,d)}(\<i,j,k\>)&=
  \left\{
  \begin{aligned}
    &w(\gamma(\xi,k),i,j),\\
    &\<-1\>,
  \end{aligned}
  \right.
&&\begin{aligned}
  &\text{if $k\le d$,}\\
  &\text{if $k>d$,}
\end{aligned}\\
  \label{eq:29}
  t_q^{\gamma(\xi,d)}(\<i,j,k\>)&=
  \left\{
  \begin{aligned}
    &w(\gamma(\xi,k),i,j\xor2^{i-1}),\\
    &\<-1\>,
  \end{aligned}
  \right.
\espc&&\begin{aligned}
    &\text{if $k\le d$,}\\
    &\text{if $k>d$,}
  \end{aligned}\\
  \label{eq:32}
  s_q^{\gamma(\eta,d)}(\<i,j,k\>)&=
  \left\{
  \begin{aligned}
    &w(\gamma(\eta,k),i,j),\\
    &\<-1\>,
  \end{aligned}
  \right.
\espc&&\begin{aligned}
  &\text{if $k\le d$,}\\
  &\text{if $k > d$,}
\end{aligned}\\
  \label{eq:33}
  t_q^{\gamma(\eta,d)}(\<i,j,k\>)&=
  \left\{
  \begin{aligned}
    &w(\gamma(\eta,k),i,j\xor2^{i-1}),\\
    &\<-1\>,
  \end{aligned}
  \right.
\espc&&\begin{aligned}
    &\text{if $k\le d$,}\\
    &\text{if $k>d$;}
  \end{aligned}
\end{alignat}
and for all $j<2^{n_q}$ and $k<n_q$: for all $d<l$, $s^{\gamma(\xi,d)}_q(\<n_q,j,k\>)$ and
$t^{\gamma(\xi,d)}_q(\<n_q,j,k\>)$ are as in~\eqref{eq:26} and~\eqref{eq:27} with
$i=n_q$, and for all $d=l,\dots,l+l^*-1$,
\begin{alignat}{2}
  s_q^{\gamma(\xi,d)}(\<n_q,j,k\>)&=
  \left\{
  \begin{aligned}
    &w(\gamma(\xi,k),n_q,j),\\
    &\<-1\>,\\
    &w(\gamma(\xi,{k-l^*}),n_q,j),\\
    &\<-1\>,
  \end{aligned}
  \right.
\espc&&\begin{aligned}
  &\text{if $k\le d$,}\\
  &\text{if $d<k<l+l^*$,}\\
  &\text{if $l+l^*\le k\le d+l^*$,}\\
  &\text{if $k> d+l^*$,}
\end{aligned}
\intertext{and $t^{\gamma(\xi,d)}_q(\<n_q,j,k\>)$ is as in~\eqref{eq:29} (with $i=n_q$), 
and $s_q^{\gamma(\eta,d)}(\<n_q,j,k\>)$ is as in~\eqref{eq:32} and}
  t_q^{\gamma(\eta,d)}(\<n_q,j,k\>)&=
  \left\{
  \begin{aligned}
    &w(\gamma(\eta,k),n_q,j\xor2^{n_q-1}),\\
    &\<-1\>,\\
    &w(\gamma(\eta,{k-l^*}),n_q,j),\\
    &\<-1\>,
  \end{aligned}
  \right.
\espc&&\begin{aligned}
    &\text{if $k\le d$,}\\
    &\text{if $d<k<l+l^*$,}\\
    &\text{if $l+l^*\le k\le d+l^*$,}\\
    &\text{if $k> d+l^*$.}
  \end{aligned}
\end{alignat}

Let us first explain why $q=(\vomega_{p(\xi)}\cup\vomega_{p(\eta)},
\{s^\alpha_q,t^\alpha_q:\alpha\in\vomega_{p(\xi)}\cup\vomega_{p(\eta)}\})$ is a
member of $\Q$. Condition~\eqref{item:1} is obvious,~\eqref{item:2} is clear, and
condition~\eqref{item:6} is satisfied since $n_q\ge l+2l^*=|\vomega_q|$. 
Condition~\eqref{item:3} holds because
$\mu([w(\alpha,i,j)])=2^{-i}$. For all $i=n+1,\dots,n_q$, for all $j<2^i$:  $d<i$
for all $d=0,\dots,l+l^*-1$ by condition~\eqref{item:6} for $p(\xi)$, and 
$s_q^{\gamma(\xi,d)}(\<i,j,d\>)=w(\gamma(\xi,d),i,j)$ and 
$t_q^{\gamma(\xi,d)}(\<i,j,d\>)=w(\gamma(\xi,d),i,j\xor2^{i-1})$
witnessing condition~\eqref{item:4} for $\gamma(\xi,d)$; similarly for $\gamma(\eta,d)$.

Observe that for all $\alpha,\beta\in\vomega_{p(\xi)}$, for all $i=n+1,\dots,n_q$, 
either $s^\alpha(\<i,j,k\>)=w(\gamma(\xi,k),i,j)$ and
$t^{\beta}(\<i,j,k\>)=w(\gamma(\xi,k),i,j\xor 2^{i-1})$, or else at least one of\linebreak
$s^\alpha(\<i,j,k\>)$ and $t^\beta(\<i,j,k\>)$ is $\<-1\>$; and similarly for all
$\alpha,\beta\in\vomega_{p(\eta)}$. First of all, taking $\alpha=\beta$ this shows
that condition~\eqref{item:7} is satisfied by $q$. It also shows that for all $\alpha,\beta\in\vomega_{p(\xi)}$, 
for all $i=n+1,\dots,n_q$ and all $j,k$,
$[s^\alpha_q(\<i,j,k\>)]\cdot[t^\beta_q(\<i,j,k\>)]=0$; 
and similarly for all $\alpha,\beta\in\vomega_{p(\eta)}$; we also see that for all
$\alpha\in\vomega_{p(\xi)}$ and $\beta\in\vomega_{p(\eta)}$, either 
$[s^\alpha_q(\<i,j,k\>)]\cdot[t^\beta_q(\<i,j,k\>)]$ is of the form
$[w(\zeta,i,j)]\cdot[w(\gamma,i,j')]$ for some $\zeta\ne\gamma$ in $\Lim(\oone)$,
which has measure $2^{-2i}$ by stochastic independence, or  it is equal to $0$;
similarly for $\alpha\in\vomega_{p(\eta)}$ and $\beta\in\vomega_{p(\xi)}$. Hence
to verify~\eqref{item:9} it remains to consider the pairs
$(\gamma(\xi,d),\gamma(\eta,\bar d))$ and
$(\gamma(\eta,d),\gamma(\xi,\bar d))$ for $d,\bar d=l,\dots,l+l^*-1$ and $i=1,\dots,n$. 
If $\delta\bigl(s^{\gamma(\xi,d)}_q(\<i,j,k\>)\bigr)
=\delta\bigl(t^{\gamma(\eta,\bar d)}_q(\<i,j,k\>)\bigr)$ then $\<i,j,k\>\in G_d$
by~\eqref{eq:23}, and hence by~\eqref{eq:24},
$s^{\gamma(\xi,d)}_q(\<i,j,k\>)=u(d,\<i,j,k\>)=s^{\gamma(\eta,d)}_q(\<i,j,k\>)$ and
thus $\mu\bigl(\bigl[s^{\gamma(\xi,d)}_q(\<i,j,k\>)\bigr]\cdot\bigl[t^{\gamma(\eta,d^*)}_q(\<i,j,k\>)\bigr]\bigr)
\le2^{-2i}$
by condition~\eqref{item:9} for $p(\eta)$; otherwise, when
$\delta\bigl(s^{\gamma(\xi,d)}_q(\<i,j,k\>)\bigr)\ne\delta\bigl(t^{\gamma(\eta,d^*)}_q(\<i,j,k\>)\bigr)$
condition~\eqref{item:9} holds by stochastic independence
and the condition~\eqref{item:3} for $p(\xi)$ and $p(\eta)$.

For condition~\eqref{item:8} we only need to look at
$\{\alpha,\beta\}\nsubseteq\vomega_{p(\xi)},\vomega_{p(\eta)}$; hence, fixing $d,\bar
d=l,l+1,\dots,l+l^*-1$ we need to consider $\gamma(\xi,d)$ and $\gamma(\eta,\bar d)$.
First suppose $d\ne\bar d$. Note that $D_{d\bar d}\subseteq F_d$ by
condition~\eqref{item:10}. Let $\Phi$ be the automorphism on $\random_{(\oone)}$
induced by $\varphi_{\gamma(\xi,d)\gamma(\eta,d)}$. Then by
conditions~\eqref{item:10},~\eqref{item:11} and~\eqref{eq:25},  
\begin{equation}
  \label{eq:37}
  [s_q^{\gamma(\xi,d)}(m)]\cdot[t_q^{\gamma(\eta,\bar d)}(m)]
  =\Phi\bigl([s_{p(\eta)}^{\gamma(\eta,d)}(m)]\cdot[t_{p(\eta)}^{\gamma(\eta,\bar d)}(m)]\bigr)
  \espc\text{for all $m\in D_{d\bar d}$}.
\end{equation}
Thus  $\sum_{m\in D_{d\bar d}}[s_q^{\gamma(\xi,d)}(m)]\cdot[t_q^{\gamma(\eta,\bar d)}(m)]
  =\Phi\bigl(\sum_{m\in D_{d\bar d}}[s_{p(\eta)}^{\gamma(\eta,d)}(m)]
    \cdot[t_{p(\eta)}^{\gamma(\eta,\bar d)}]\bigr)
  =\Phi(1)=1$, as required. Otherwise, when $d=\bar d$, put  
\begin{align}
  \label{eq:21}
  \begin{split}
  H_d&=\{\<i,j,l_d(i,j)\>:i=1,\dots,n\text{, }j<2^i\}\\&\qquad
\cup\{\<i,j,d\>:i=n+1,\dots,n_q\text{, }j<2^i\},
  \end{split}\\
  I_d&=\{\<n_q,j,d+l^*\>:j<2^{n_q}\}.
\end{align}
Then $\bigl\{[s_q^{\gamma(\xi,d)}(m)]\cdot [t_q^{\gamma(\eta,d)}(m)]:m\in H_d\bigr\}
=\bigl\{[w(\gamma(\xi,d),i,j)]\cdot[w(\gamma(\eta,d),\allowbreak i,j)]:
i=1,\dots,n_q\text{, }j<2^i\bigr\}$ and 
$\bigl\{[s_q^{\gamma(\xi,d)}(m)]\cdot[t_q^{\gamma(\eta,d)}(m)]:m\in I_d\bigr\}
=\bigl\{[w(\gamma(\xi,d),n_q,j)]\cdot[w(\gamma(\eta,d),n_q,j)]:j<2^{n_q}\bigr\}$, and
therefore $H_d\cup I_d$ satisfies~\eqref{item:8} for $q$ by~\eqref{eq:13};
similarly, for the opposite pair. 

It should be clear now that in fact $q\le p$, concluding the proof.
\end{proof}

Now for a filter $\G\subseteq\Q$, define
\begin{equation}
  \label{eq:18}
  s^\alpha=\bigcup_{p\in\G}s^\alpha_p\And t^\alpha=\bigcup_{p\in\G}t^\alpha_p.
\end{equation}
For a sufficiently generic $\G\subseteq\Q$, a similar argument as in the proof
of Lemma~\ref{l-3} shows that for each $\alpha$,
$s^\alpha,t^\alpha:\N\to\Fin(\alpha+\omega,\two)\cup\{\<-1\>\}$. 
Define  $\random_{(\oone)}$-names 
$(\dot a_\alpha,\dot b_\alpha:\alpha<\oone)$ by
\begin{align}
  \label{eq:10}
  \lbrak m\in \dot a_\alpha\rbrak&=[s^\alpha(m)],\\
  \label{eq:15}
  \lbrak m\in \dot b_\alpha\rbrak&=[t^\alpha(m)],
\end{align}
for all $m$. Condition~\eqref{item:5} ensures that $\dot a_\alpha$ and $\dot
b_\alpha$ are increasing with respect to $\alpha$ modulo $\Fin$. By~\eqref{item:7}, 
$\lbrak\dot a_\alpha\cap\dot b_\alpha=\emptyset\rbrak=1$ for all $\alpha$ establishing~\eqref{eq:3}.
For all $\alpha\ne\beta$, by~\eqref{item:8}, 
$\lbrak\dot a_\alpha\cap\dot b_\beta\ne\emptyset\rbrak=1$
establishing~\eqref{eq:5}. And~\eqref{item:9} establishes~\eqref{eq:50}.
This completes the proof that the object described in~\eqref{eq:3}--\eqref{eq:50} is
consistent, and thus completes the proof of the Theorem.

\bibliography{Hausdorff.gap}

\providecommand{\bysame}{\leavevmode\hbox to3em{\hrulefill}\thinspace}
\providecommand{\MR}{\relax\ifhmode\unskip\space\fi MR }
\providecommand{\MRhref}[2]{%
  \href{http://www.ams.org/mathscinet-getitem?mr=#1}{#2}
}
\providecommand{\href}[2]{#2}
\begin{thebibliography}{Dow95}

\bibitem[AT97]{rst:A}
Uri Abraham and Stevo Todor{\v{c}}evi{\'c}, \emph{Partition properties of
  $\omega\sb 1$ compatible with {C}{H}}, Fund. Math. \textbf{152} (1997),
  no.~2, 165--181.

\bibitem[Dow95]{D}
Alan Dow, \emph{More set-theory for topologists}, Topology Appl. \textbf{64}
  (1995), no.~3, 243--300.

\bibitem[Fre01]{F3}
David~H. Fremlin, \emph{Measure algebras}, vol.~3, Torres Fremlin, 2001,
  \texttt{http://www.essex.ac.uk/maths/staff/fremlin/mtsales.htm}.

\bibitem[Gil36]{G}
J.~Gillis, \emph{Note on a property of mesurable sets}, J. London Math. Soc.
  \textbf{11} (1936), 139--141.

\bibitem[GS01]{GS}
Moti Gitik and Saharon Shelah, \emph{More on real-valued measurable cardinals
  and forcing with ideals}, Israel J. Math. \textbf{124} (2001), 221--242.

\bibitem[Had94]{Had}
Jacques Hadamard, \emph{Sur les caracteres de convergence des series a termes
  positifs et sur les fonctions indefiniment croissantes}, Acta Math.
  \textbf{18} (1894), 319--336.

\bibitem[Hau36]{Ha}
F.~Hausdorff, \emph{Summen $\aleph_1$ von {M}engen}, Fund. Math. \textbf{26}
  (1936), 241--255.

\bibitem[Hir01]{HI}
James Hirschorn, \emph{Summable gaps}, to appear in Annals of Pure and Applied
  Logic, 2001.

\bibitem[Kei70]{Ke}
H.~Jerome Keisler, \emph{Logic with the quantifier ``there exist uncountably
  many''}, Ann. Math. Logic \textbf{1} (1970), 1--93.

\bibitem[Kun76]{Ku}
Kenneth Kunen, \emph{$(\kappa,\lambda^*)$ gaps under $\ma$}, handwritten note,
  August 1976.

\bibitem[Sch93]{Sch}
Marion Scheepers, \emph{Gaps in $\omega\sp \omega$}, Set theory of the reals
  (Ramat Gan, 1991), Bar-Ilan Univ., Ramat Gan, 1993, pp.~439--561.

\bibitem[Sch96]{Sch1}
\bysame, \emph{The {B}oise problem book}, Topology Atlas Preprint \#167, 1996.

\bibitem[TF95]{TF}
Stevo Todor{\v c}evi\'c and Ilijas Farah, \emph{Some applications of the method
  of forcing}, Yenisei, Moscow, 1995.

\end{thebibliography}
\bibliographystyle{amsalpha}

\end{document}